# A two-dimensional introduction to sashiko

Carol Hayes[1] and Katherine A. Seaton[2]

[1]School of Culture, History and Language, CAP, ANU, ACT 2600 Australia;
carol.hayes@anu.edu.au
[2]Department of Mathematics and Statistics, La Trobe University, VIC 3086 Australia;
k.seaton@latrobe.edu.au

## Abstract

Through the hands-on creation of two sashiko pieces of work – a counted thread kogin bookmark and a single running stitched hitomezashi sampler – participants will explore not only the living cultural history of this traditional Japanese needlework but will also experience the mathematics of sashiko in a tangible form, and will take away with them items of simple beauty.

## About this workshop

Sashiko is inherently two-dimensional, being a form of stitching on fabric, or through layers of fabric. This workshop will also be two-dimensional in the sense that the presenters approach sashiko from different directions. Carol Hayes' area of expertise is Japanese language and culture, while Katherine Seaton is a mathematician. A third two-dimensional element speaks to Japanese aesthetic tradition. In the eyes of renowned Edo poet Matsuo Bashō, the vertical axis represents a link back to past tradition while the horizontal places the work in the everyday world of the "contemporary, urban commoner life and a new social order" [16]. So the vertical is an unbreakable link to the sincere heart of the cultural past and the horizontal a challenge to those conventions, offering a new contemporary perspective [4]. Contemporary sashiko is the union of these two dimensions.

This appears to be only the second appearance of sashiko at a Bridges conference. In 2006, Setsu Pickett described geometric features of contemporary pattern sashiko (*moyōzashi*), such as parallel lines, grids and arcs of circles, in a short talk [13], and illustrated them in her art exhibition piece *Counterweights and Plovers* [14]. Her motivation to analyse sashiko, and the underlying geometry of the hemp leaf pattern (*asa-no-ha*) in particular, was to use it as an inspiration for her usual medium of velvet weaving designs [15]. In this current workshop, participants will make their own small pieces of sashiko.

## Sashiko

Sashiko needlework began as a functional, simple running stitch used to repair or strengthen garments, to patch worn clothes and to quilt multiple layers of rough fabric together for warmth. Translating as 'little stabs' the name refers to the act of pushing the needle through cloth. Sashiko has its origins in the practical rather than the decorative, born of the "wisdom of everyday life" embodying the "beauty of the functional" [10]. The simple and repetitive technique of sashiko stitching requires the stitcher to be present in the moment, to focus on the work at hand.

The tradition of patching fragments of material together to create a durable fabric, which is to be found around the world, is a means of salvaging a scarce resource. Rags were collected and stored from older items, and scraps were bought from merchants selling used garments from the more affluent Edo, Osaka and Kyoto regions [20]. This method of darning and quilting with running stitches is also referred to as *boro* and it has become increasingly well-known in recent years in the slow fashion movement as a form of visible mending [5]. See Figure 1 for a traditional example.

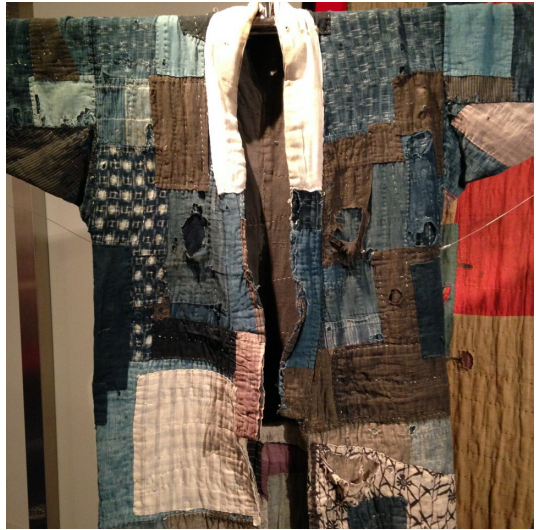

**Figure 1:** *Boro kimono, Amuse Museum Permanent Exhibition. (Photo by Carol Hayes.)*

The trajectory of sashiko is like that of stitching across the world. As embroidery researcher Lanto Synge notes: "The earliest needlework was of a plain, practical nature, done with strong fibrous materials such as hair, to join skins and furs for clothing, and embroidery was used to strengthen parts subject to greater wear. From this basic necessity a sumptuous decorative art gradually emerged" [19].

There are three main types of traditional sashiko; *hishizashi*, *koginzashi*, and *shōnai sashiko*. The dense quilting stitches of Tsugaru Koginzashi and Nambu Hishizashi both originated in Aomori prefecture in the Northern Tohoku region of Japan where the weather could be very cold, and thick, sturdy fabric was needed for clothes. In the shōnai form *hitomezashi*, which means "one-stitch", the pattern is created by using relatively large stitches. These styles are shown in Figure 2.

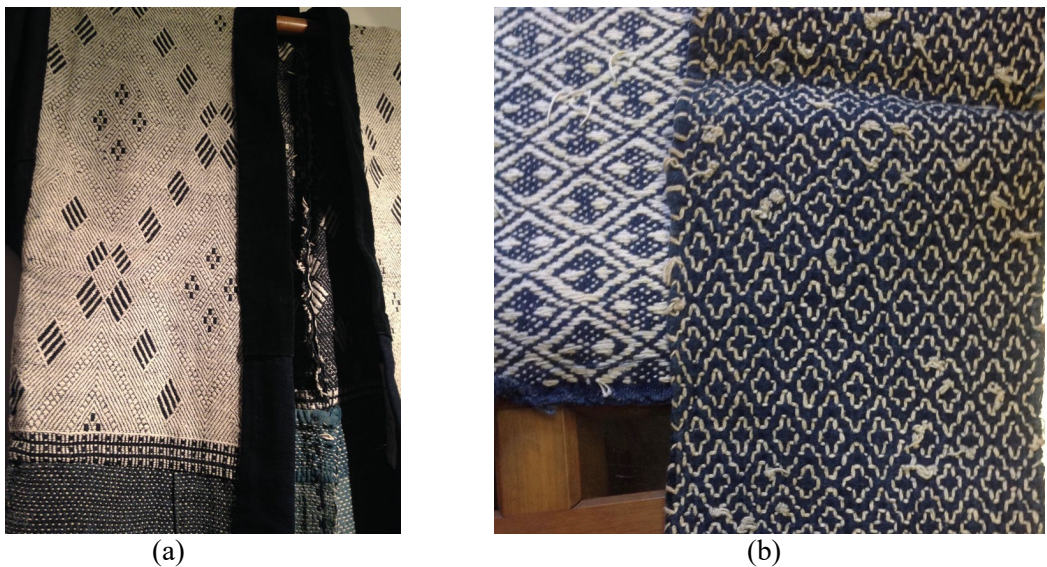

(a) (b)

**Figure 2:** *Traditional sashiko: (a) Nambu hishizashi, Amuse Museum collection, (b) Antique Tsugaru koginzashi, behind a piece of shōnai hitomezashi, Yoshiura Collection. (Photos by Carol Hayes.)*

All three styles employ counted-thread techniques wherein an underlying grid, provided by the warp and weft of the fabric, is used to create the design, as in cross-stitch or black work. Unlike cross-stitch, however, which is generally worked from the centre of the piece outwards, sashiko is worked from edge to edge. One of the main differences between the kogin and hishi styles is that the kogin patterns use an odd number count of fabric threads, while the hishi pattern is based on even thread count. In some hitomezashi patterns, stitches are worked one at a time and meet at points on the underlying grid; in the other styles the space between stitches remains visible.

Contemporary pattern sashiko (*moyōzashi*) pieces use spaced running stitches on a single layer of fabric, pre-printed or marked up using a template prior to stitching. This form of sashiko has been popularised globally, see for example [3]. The hitomezashi kit guide [9] in Figure 3 demonstrates that, after working the border, the vertical (1), horizontal (2) and diagonal (3) lines are completed in ordered progression using uniform stitches.

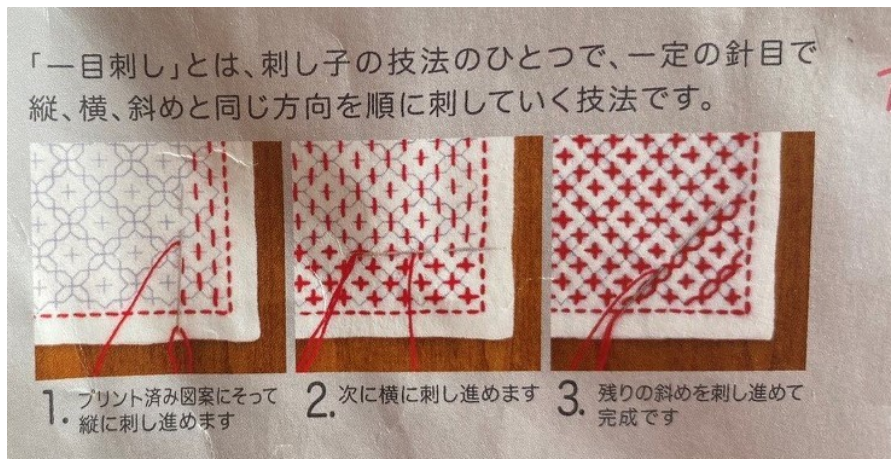

**Figure 3:** *A modern commercially produced sashiko kit [9], showing both the fabric pre-marked with the stitch placement, and the instruction to work all stitches in a given direction before moving to the next direction. (Photo by Carol Hayes.)*

## Sashiko and Japanese Culture

Hierarchy was an important element of pre-modern societies, and Japan was no different. The story of sashiko is also a story of self-worth, of how one accounts for oneself. Clothing, both functional and decorative, speaks to both the status and the taste of the wearer. Sashiko is also a narrative of the human desire for beauty, for the decorative in the most impoverished functional space [8].

Although running stitch has a long history in Japan (as in the rest of the world), the distinctive sashiko stitching patterns developed in the Edo period (1603-1868). During this period, which saw the rise of the merchant middle class, the government released a series of increasingly punitive sumptuary laws and regulations limiting private expenditure on clothing and other personal items. The aim was to regulate differences between social classes by curbing conspicuous consumption, particularly when it caused blurring of the lines between social classes [12].

These laws provide an important context for the development of sashiko stitching, because under the this regime both silk and cotton were forbidden to the commoners [17]. The people of the Tohoku region developed their own cloth using hemp, usually dyed with indigo, which is rough and has an open weave. It can be argued that surface sashiko stitching developed to soften the poor quality of the fabric used, and the quilting allowed wadding from plant matter or old rags to be added between layers for further warmth. Certain occupations caused different types of wear-and-tear to garments influencing pattern placement; farmers and labourers carried heavy loads on their shoulders, while firemen needed thick coats to protect them from the flames. Patching was originally done with thread that came from the garment itself, but as undyed white cotton thread became available it began to be used to mend garments. Although the white thread faded with use, each layer of sashiko mending allowed for the potential of patterns to develop in the stitching. The running stitch lent itself to stop-start working, and so commoners who had saved even small lengths of thread could use this style of patterned stitching.

As the practice of functional sashiko stitching became more widespread, decorative patterns developed and it began to be used not only by the peasant classes but by lower and middle classes as a more decorative form of domestic needlework. Household wash cloths and hand towels began to be stitched with sashiko patterns. Sashiko developed into a decorative stitching that uses numerous patterns to decorate clothes and daily items such as napkins. These patterns reflect the nature of Japanese society: farming and fishing, the seasons and the natural landscape of Japan. Traditionally, and now as part of the mindfulness movement, sashiko stitching can be a meditative experience of well-being, with each stitch a deliberate thought or small prayer [8].

## Mathematics of Sashiko

There are obvious decorative geometric patterns in the pieces shown in Figure 2, or in Figure 4, but this is just the surface mathematics of sashiko. We will limit our attention to hitomezashi on the square grid and find much rich content even in this confined environment.

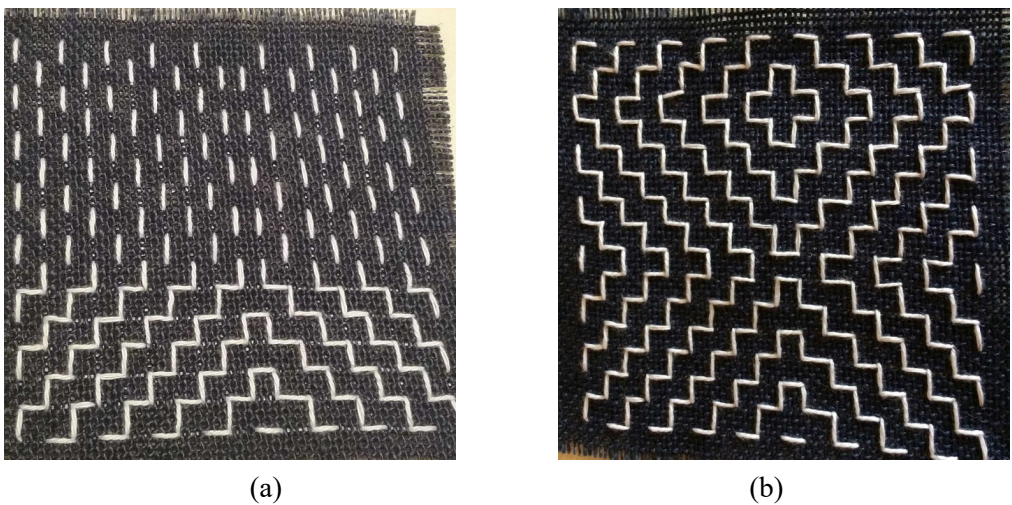

(a) (b)

**Figure 4:** *A hitomezashi sampler as will be made in this workshop: (a) horizontal stitches being added to the vertical stitches already in place, revealing the pattern; (b) mountain form (yamagata) in the bottom of the piece, and a variant of the persimmon flower (kakinohanazashi) towards the top.*
*(Work and photo by K. Seaton.)*

### *Counting*

Consider the lines of running stitch, shown in Figure 4(a). There are only two ways in which the first stitch in a line can be made; either it is formed on the front of the fabric, or on the reverse. There is then, by the nature of running stitch, no choice about how to place the subsequent stitches in that line. Each line can be in only one of two (binary) states. Thus, if $m$ vertical lines and $n$ horizontal lines comprise a piece of hitomezashi, the number of ways in which it could be formed is $2^{m+n}$. Many of these 'designs', however, would be unrecognizable as traditional sashiko!

### *Duality*

Duality is a deep and pervasive concept in modern mathematics, although its realisation can be very different depending on context. In sashiko, duality manifests itself in the complementary patterns that are formed by the stitches on the two sides of the stitched fabric. In Figure 5, it can be seen that the pattern of crosses (which resemble the kanji character for the number ten 十) has as its dual a pattern of squares and zig-zagging stepped lines.

*Symmetry*

Hitomezashi worked on the square grid is a constrained fibre art form. Unlike cross-stich, in which one of the two stitches making up the cross necessarily lies above the other, so that symmetry is broken if examined at the stitch level, in hitomezashi the stitches lie flat on the fabric. Thus this embroidery form has the potential to give a beautiful, more perfect, realisation of the rosette, frieze and wallpaper symmetry groups studied previously in the context of other fibre arts [7]. In Figure 4(b) we can see a vertical axis of symmetry (only), whereas in Figure 2(b) there is all the symmetry of the wallpaper group cmm. Designs with 90° rotations are also possible. It is perhaps emblematic of the differing mindsets that a mathematician and a scholar of culture bring to sashiko that Katherine immediately remarked upon the four-fold rotations of many designs, whereas Carol, knowing that even numbers are considered less favourably than odd numbers in Japanese culture, notes that the crosses contain five squares arranged in three rows as 1-3-1.

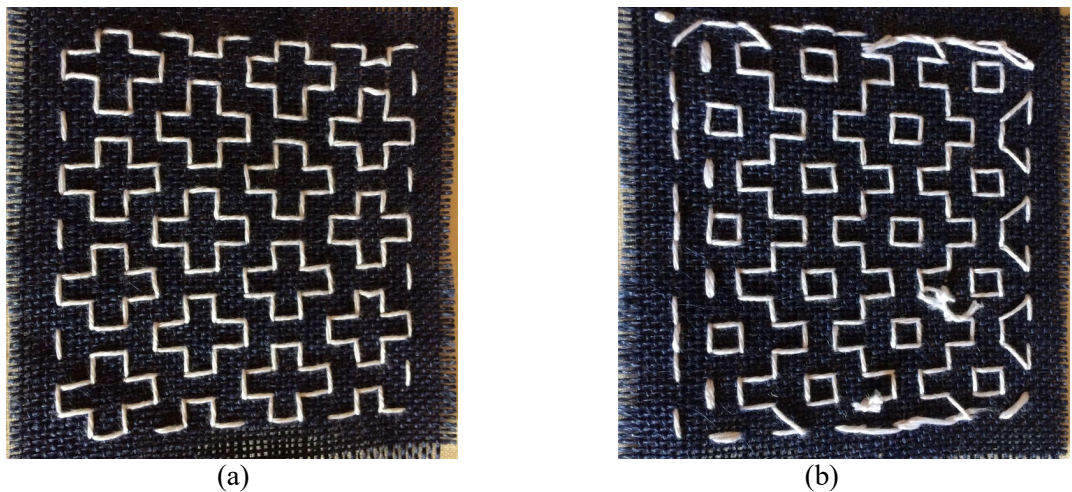

(a)          (b)

**Figure 5:** *Duality: (a) 'ten' crosses (jūjizashi) on one side of the work, and (b) the dual design of squares and stepped lines on the reverse. (Work and photo by K. Seaton.)*

*Loops, Truchet tiles, snowflakes and polyominoes*

When there is no overall wallpaper-type structure to a piece of hitomezashi, as in the piece shown in Figure 6, there can be both nested loops (polyominoes) and stepped lines of stitches that traverse the fabric from edge to edge. Such structures are found in the fully-packed loop models of statistical mechanics, with connections to alternating sign matrices, percolation and the six-vertex model [6]. They also resemble Truchet tilings [18]; their classification lies in the field of enumerative combinatorics.

Remarkably, the loops or polyominoes (that is, objects formed by the union of edge-adjacent lattice squares) of hitomezashi have a connection to Fibonacci. The so-called Fibonacci snowflakes (which do not have the six-fold structure of a true snowflake, but which have fractal features in common with the better-known Koch snowflake, that does) arise by generating a path of turns on the square lattice using the Fibonacci word [1]. These snowflakes, and their generalisations [11], tile the plane [2]. Among the Fibonacci snowflakes of Figure 7, see the square, the ten-cross, the outline of the persimmon, and another design that can be realised using superposed horizontal and vertical running stitch. The number of stitches in the outline of a Fibonacci snowflake is an odd Fibonacci number (1, 3, 5, 13) multiplied by four.

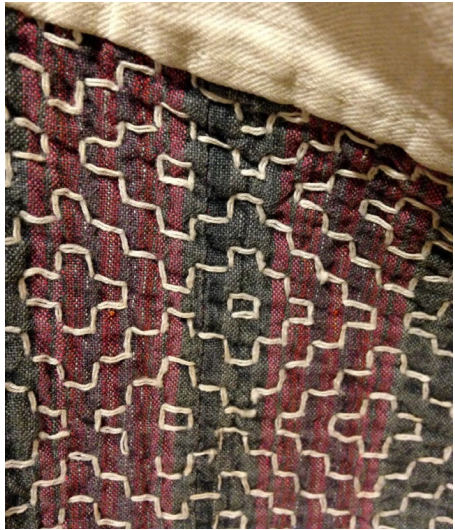

**Figure 6:** *Hitomezashi apron detail, Amuse Museum collection. (Photo by Carol Hayes.)*

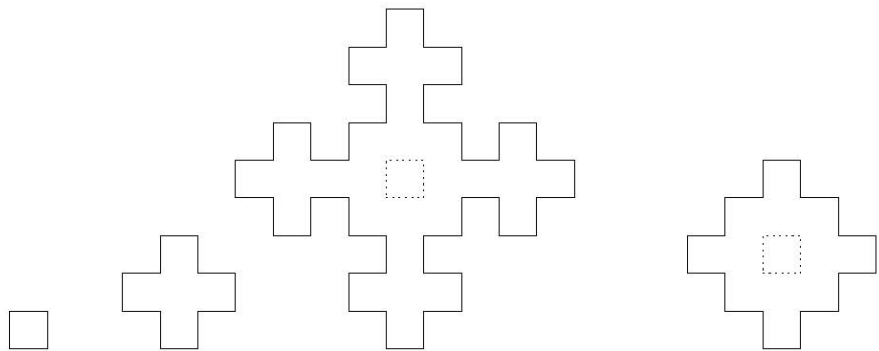

**Figure 7:** *Fibonacci snowflakes [1] of order 0 (square),1(jūjizashi) and 2, and to the right, a generalised snowflake [11]. The dashed lines indicate how these structures arise with a nested square when stitched using running stitch in hitomezashi; the figure on the right is the traditional persimmon flower design.*

## The Activity

The core of this workshop is the following activity that engages participants in both the cultural practice and mathematical aspects of sashiko. They will be provided with instructions and materials to produce two small pieces of sashiko (or a drawing of each), one a sampler in the hitomezashi style, which would be suitable for a coaster, and the other a bookmark in the kogin style.

Like many fibre arts projects, simple sashiko can be stitched while listening to music or conversation. Once the activity is underway, Carol will provide an account of how a craft born of frugality has survived mechanisation of clothing production to become a decorative and valued art form. Katherine will facilitate discussion of mathematical features of hitomezashi and kogin sashiko beyond the obvious geometric ones, drawing on the mathematical expertise among participants.

Each participant will be supplied with a round-ended (tapestry) needle, white cotton thread, printed patterns and instructions, and a 2.5 x 8 cm rectangle of coloured hessian (for the bookmark) and a 12 cm square of blue hessian (for the coaster). The choice of hessian is quite deliberate, both for practicality (it being relatively inexpensive and having a regular open weave) and for reasons that will become apparent as Carol speaks.

*Hitomezashi Coaster*

The first page of the pattern will show the vertical running stitches to be worked, and the second page the horizontal stitches. Running stitch (*yokogushi*) is a common sashiko stitch pattern in its own right. An element of wonder and surprise is that the two sets of stitches will not be provided as a superimposed image, as they will appear in the final piece. In order to discover the design they have been given to create, participants will have to complete both sets of stitches. There will be a variety of different patterns distributed across the group. Two or more may have the same vertical stitch pattern, but a different horizontal stitch pattern, and vice versa.

*Kogin Bookmark*

Participants will choose from four traditional kogin patterns – dragonfly, butterfly, gourd and *kikurako* (sitting crossed legged on the floor), as shown in Figure 8. Working only on the horizontal, stitching will be completed in stitch lengths of 1, 3 and 5 to recreate the template on a small hessian rectangular bookmark.

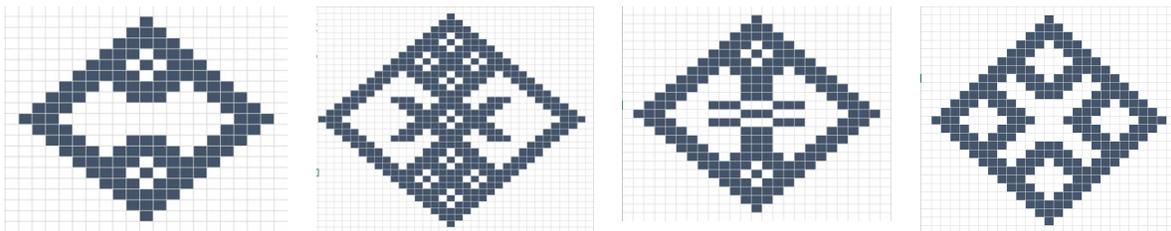

**Figure 8:** *Kogin motifs for the bookmarks to be made in the workshop, which represent natural elements worked with horizontal stitches of odd length. From left: gourd, butterfly, dragonfly and crossed-legs.*

In each case, the written instructions will be supplemented with verbal instructions and hands-on assistance. If there are participants who do not wish to sew for any reason, graph paper, pens and rulers will also be available for them to transcribe patterns onto a diagram, to experience the emergence of the overall design. By the end of the workshop, each participant will have produced two small pieces of sashiko (or two drawings). We will reveal which participants had the same vertical or horizontal stitches within their hitomezashi patterns, and will encourage comparison of the final pieces.

## Summary and Conclusions

In this workshop, participants will have created their own pieces of sashiko, and will have learned about its history, cultural and artistic properties, and its mathematical features. There is more to be said about this last aspect, not only for hitomezashi, but also for other traditional, constrained variants of sashiko, such as kogin and hishizashi. Indeed, from their own mathematical backgrounds, workshop participants may be able to identify other mathematical connections.

We hope the participants will enjoy the workshop activity on a number of different levels: as a part of a living cultural tradition with roots deep in Japanese daily tradition, and as a product of both beauty and mathematical complexity.

## Acknowledgements

KAS acknowledges helpful [e]-conversations with Jan de Gier and Bernard Nienhuis (mathematical physics) and Julia Collins (Truchet tilings). Michael Assis (origami designer and physicist) made a remark about Aztec diamonds which led her to Fibonacci snowflakes.